\documentclass{jnmp}

\usepackage{amsmath}

\JNMPnumberwithin{equation}{section}

\theoremstyle{definition}

\newtheorem*{remark}{Remark}

\begin{document}

%
\renewcommand{\evenhead}{I~M~Nefedov and I~A~Shereshevskii}
\renewcommand{\oddhead}{Operator Exponential Method for Initial-Boundary Value
Problems}

%
\thispagestyle{empty}

\FirstPageHead{8}{3}{2001}{\pageref{nefedov-firstpage}--\pageref{nefedov-lastpage}}{Letter}

\copyrightnote{2001}{I~M~Nefedov and I~A~Shereshevskii}

\Name{Solving the Difference Initial-Boundary Value Problems
by the Operator Exponential Method}
\label{nefedov-firstpage}

\Author{I~M~NEFEDOV and I~A~SHERESHEVSKII}

\Address{Institute for Physics of Microstructures, RAS,\\
GSP-105, 603950 Nyzhny Novgorod, Russia\\
E-mail: ilya@ipm.sci-nnov.ru; nefd@ipm.sci-nnov.ru}

\Date{Received September 19, 2000; Revised February 15, 2001;
Accepted March 22, 2001}

\begin{abstract}
\noindent
We suggest a modification of the operator exponential method for the
numerical solving the difference linear initial boundary value
problems.  The scheme is based on the representation of the difference
operator for given boundary conditions as the perturbation of the same
operator for periodic ones.  We analyze the error, stability and
efficiency of the scheme for a model example of the one-dimensional
operator of second difference.
\end{abstract}

\section{Introduction}

Numerical solution of the linear difference initial-boundary value
problems is an essential part of modelling of the physical processes
and phenomena described by the evolutionary differential equations,
such as the Schr\"odinger equation,  diffusion equation,  Ginzburg--Landau
equation and many other.

Along with the classical grid methods~\cite{nefedov:1},  an ever increasing use
in treatment of such evolutionary problems is currently made of the
operator exponential (OE) method~\cite{nefedov:2},  which is based on the
Lee--Trotter--Kato formula~\cite{nefedov:3} for approximate calculation of the
exponential of the sum of noncommuting matrices.

The OE method offers a number of advantages relevant to both explicit
and implicit difference schemes.  It does not involve iteration
procedures and often proves to be absolutely stable.  Its
applicability is, however, limited by the impossibility to explicitly
calculate the exponential of the difference operators expressed in a
general form.  In fact, effective algorithms of exponential
calculation exist only for the difference operators with constant
coefficients and periodic boundary conditions on ``rectangular"
subsets of ${\mathbb Z}^n$.  These algorithms are based on the fast Fourier
transform~\cite{nefedov:4} and allow one to calculate the exponential in
$O(N\log_2N)$ operations, where $N$ is the number of points in the
domain.  For other boundary conditions such algorithms are not
available.

In this work a linear difference operator with assigned boundary
conditions is considered as perturbation of the same operator
with periodic boundary conditions,  and the exponential of
such an operator is calculated by the Lee--Trotter--Kato formula.

The perturbating operator is essentially an operator in the space of
functions on the domain's boundary, so the problem of calculating its
exponential is essentially simpler because the number of boundary
points is generally much smaller than the total number of points in
the domain.  This ensures practically the same efficiency of the
algorithm proposed as that obtained with the Fourier method for
solving the periodic boundary value problems.

Analysis of the error, stability and efficiency of the algorithm
proposed is generally quite complicated.  So we only present it for a
model example of the operator of second difference (one-dimensional
difference Laplace operator).  This case is probably least
``favorable" for the OE method due to availability of effective
difference schemes such as the {\it sweep method}~\cite{nefedov:1}.
Nevertheless, the algorithm proposed is competitive with the
well-known schemes, in particular, as applied to the Schr\"odinger
equation.

The idea of representing differential operators with various boundary
conditions as one another's perturbations was put forward by
M~G~Krein~\cite{nefedov:6} and is being actively used in modern mathematical
physics (see, for example,~\cite{nefedov:7}).  Applicability of the Krein
method to difference operators was considered in~\cite{nefedov:8}.
The results of this work were reported at the ``Conference on
differential equations and applications" (Saransk, Russia, 1994).  The
summary of this report was published in \cite{5}.

\section{Description of the method}

Let ${\mathcal M}$ be a set of points and $C({\mathcal M})$ a set of
complex-valued functions on $\mathcal M$.  Let $\hat{A}:C({\mathcal
M})\longrightarrow C({\mathcal M})$ be a linear operator of the form
\begin{equation}
(\hat{A}f)(x)=\sum_{y\in\gamma_x}a_x(y)f(y), \qquad f\in C({\mathcal M}),
\label{nefedov:2.1}
\end{equation}
where $\gamma_x$ is a finite subset of $\mathcal M$ for each value of
$x$, and $a_x(\cdot)$ is a given function on $\gamma_x$.

Let $\Omega$ be a subset of $\mathcal M$.  We call point $x\in\Omega$
an {\it inner} point of $\Omega$ {\it relative to} $\hat{A}$ if
$\gamma_x\subseteq\Omega$, and a {\it boundary} point of $\Omega$ {\it
relative to} $\hat{A}$ if $\gamma_x$ does not completely lie in
$\Omega$.  Denote by $\partial_A\Omega$ the set of all boundary points
of $\Omega$ relative to $\hat{A}$ and let
$b_A\Omega=\mathop{\cup}\limits_{x\in\partial_A\Omega}\gamma_x\backslash\Omega$.

Note that by definition (\ref{nefedov:2.1}), to calculate the values
of $\hat{A}f$ at boundary points of~$\Omega$, we have to know the
values of function $f$ on the set $\Omega\cup b_A\Omega$.  A linear
operator 
$$
\hat{L}: C(\Omega)\longrightarrow C\left(\Omega\cup
b_A\Omega\right)\text{ such that $(\hat{L}f)(x)=f(x)$ for all $x\in\Omega$}
$$
will be called an {\it extension operator for} $\hat{A}$.  An operator
$\hat{A}_L: C(\Omega)\longrightarrow C(\Omega)$ such that
$(\hat{A}_Lf)(x)=(\hat{A}\hat{L}f)(x)$ for all $x\in\Omega$ will be
called an {\it $L$-expansion of operator} $\hat{A}$.  The operator
$\hat{L}$ plays the same role for difference operators as the boundary
conditions play for differential operators.

We now consider a difference initial-boundary value problem for
operator $\hat{A}$:
\begin{equation}
\left\{  \begin{array}{l}
             \displaystyle \frac{\partial f}{\partial t}
                                             =\hat{A}_Lf,  \quad t\geq 0,   \vspace{2mm}\\
              f(0, x)=g(x),
         \end{array}
\right.
\label{nefedov:2.2}
\end{equation}
where $f(t, \cdot), g\in C(\Omega)$,  $\hat{L}$ is a given extension
operator for  $\hat{A}$. The solution of problem (\ref{nefedov:2.2})
is of the form:
\begin{equation}
f(t, \cdot)=\exp(t\hat{A}_L)g,
\label{nefedov:2.3}
\end{equation}
where operator $\exp(t\hat{A}_L)$ can be defined as a matrix power
series since $\Omega$ is finite.  Given~$\hat{A}$ and $\Omega$, the
efficiency of computation of $\exp(t\hat{A}_L)$ in (\ref{nefedov:2.3})
may largely depend on the extension operator $\hat{L}$.  Let us
clarify the above said with a test example.

Let $\mathcal M={\mathbb Z}$,  let
$\hat{\Delta}$ be the difference Laplacian~\cite{nefedov:1}
defined by the relation
\begin{equation}
(\hat{\Delta}f)(x)=f(x-1)-2f(x)+f(x+1),  \qquad x\in{\mathbb Z}.
\label{nefedov:2.4}
\end{equation}
In this case $\gamma_x$ in (\ref{nefedov:2.1}) is the set $\{x-1, x, x+1\}$,
\[
a_x(y)=\left\{
\begin{array}{rll}
1 &  \text{if} & y=x-1, \\
-2 &  \text{if} &  y=x,\\
1 &  \text{if} &  y=x+1.
\end{array}
\right.
\]
Let $\Omega=\{0, 1, \ldots, N-1\}$.  Then
$\partial_{\Delta}\Omega=\{0, N-1\}$, and $b_{\Delta}\Omega=\{-1,
N\}$.  Let the extension operator $\hat{L}$ correspond to the periodic
boundary conditions for $\hat{\Delta}$:
\begin{equation}
(\hat{L}f)(x)=
\left\{
\begin{array}{lll}
f(N-1) & \text{if} & x=-1,\\
f(x) & \text{if}  & x\in\Omega,\\
f(0) & \text{if}  & x=N.
\end{array}
\right.
\label{nefedov:2.5}
\end{equation}
The exponential $\exp(t\hat{\Delta}_L)$ in (\ref{nefedov:2.3}) can be expressed
by the following formula:
\begin{equation}
\exp(t\hat{\Delta}_L)=\hat{F}^{-1}\exp(t\hat{\Lambda})\hat{F},
\label{nefedov:2.6}
\end{equation}
where $\hat{\Lambda}$ is the diagonal operator of the form:
\[
(\hat{\Lambda}f)(x)=\nu(x)f(x), \qquad \nu(x)=-4\sin^2\frac{\pi x}{N}, \qquad x\in\Omega,
\]
$\hat{F}$ is the operator of the discrete Fourier transform:
\[
(\hat{F}f)(x)=\sum_{y\in\Omega}\exp\left(-\frac{i 2\pi xy}{N}\right)f(y).
\]

Note that computation of the vector $\exp(t\hat{\Delta}_L)f$ via
formula (\ref{nefedov:2.6}) takes about $N\log_2N$ operations if we make use
of the known Fast Fourier Transform (FFT) algorithm~\cite{nefedov:4}.

Let $\hat{K}$ be the extension operator for $\hat{\Delta}$,
corresponding to the boundary conditions of the 3rd kind,  i.e.,
\begin{equation}
(\hat{K}f)(x)=
\left\{
\begin{array}{lll}
\alpha f(0) & \text{if} & x=-1,\\
f(x) & \text{if} & x\in\Omega,\\
\beta f(N-1) & \text{if} & x=N,
\end{array}
\right.
\label{nefedov:2.7}
\end{equation}
where $\alpha$ and $\beta$ are,  generally,  the complex coefficients.
(The case $\alpha=\beta=-1$ corresponds to the Dirichlet boundary
conditions,  and $\alpha=\beta=1$ to the Neumann boundary conditions.)
In this case the known algorithms for exact computation of the vector
$\exp(t\hat{\Delta}_K)f$ (for example,  using expansion in
eigenfunctions of $\hat{\Delta}_K$) involve $\sim N^2$ operations.

Considering the general case again,  the question arises: whether the
available effective algorithm for the $\exp(t\hat{A}_L)$ computation
($\hat{L}$ is the given extension operator) can be used to
approximately evaluate $\exp(t\hat{A}_K)$ for another extension
operator $\hat{K}$?

Below we describe a version of an OE method which establishes the
relation between the exponents of different extensions of a difference
operator and thus answer the above question.

Let $\hat{K}$ and $\hat{L}$ be two different extension operators for
the operator $\hat{A}$.  We further assume for simplicity that these
operators satisfy the following additional condition: the equations
$\hat{L}f=\hat{L}g$ and $\hat{K}f=\hat{K}g$ are fulfilled for any $f,
g\in C(\Omega)$ such that $f(x)=g(x)$ for
$x\in\Omega\backslash\partial_A\Omega$.

Consider operator $\hat{G}_{KL}=\hat{A}_K-\hat{A}_L$.  It follows from
definition of extension operators, that $(\hat{G}_{KL}f)(x)=0$ at all inner
points $x\in \Omega$, and that $\hat{G}_{KL}f=\hat{G}_{KL}g$ if $f(x)=g(x)$ at
the inner points $x\in\Omega$.  Therefore, $\hat{G}_{KL}$ is the direct sum of
the zero operator in the subspace $C(\Omega\backslash\partial_A\Omega)$ of
$C(\Omega)$ and an operator in the subspace $C(\partial_A\Omega)$; we will
denote the restriction of $\hat{G}_{KL}$ on $C(\partial_A\Omega)$ by the same
character $\hat{G}_{KL}$. 

This suggests that when the number of boundary points of $\Omega$ is
much smaller than the total number of points in $\Omega$, the problem
of computing $\exp(t\hat{G}_{KL})$ becomes much simpler than the
initial problem of evaluating $\exp(t\hat{A}_K)$.

\begin{remark}
For the extension operators of the general form a small
modification of these arguments leads to the same result.  In the
above example of operator $\hat{\Delta}$ the number of boundary points
is equal to two, and the computation reduces to finding the exponential
of a $2\times 2$ matrix.
\end{remark}

Since $\hat{A}_K=\hat{A}_L+\hat{G}_{KL}$, the following
relations hold:
\begin{gather}
\exp(t\hat{A}_K)=\exp(t\hat{A}_L)\exp(t\hat{G}_{KL})+O\left(t^2\right)
\equiv \hat{S}_1(t)+O\left(t^2\right)
\label{nefedov:2.8}\\
\exp(t\hat{A}_K)=\exp\left(\frac12t\hat{G}_{KL}\right)
                 \exp(t\hat{A}_L)
                 \exp\left(\frac12t\hat{G}_{KL}\right)+O\left(t^3\right)\nonumber\\
\phantom{\exp(t\hat{A}_K)}\equiv \hat{S}_2(t)+O\left(t^3\right).
\label{nefedov:2.81}
\end{gather}
They are similar to the conventionally used OE schemes of
the 1st and 2nd order approxi\-mation~\cite{nefedov:2}.

Owing to the above mentioned properties of operator $\hat{G}_{KL}$,
these formulas allow one to roughly calculate the exponential of
$\hat{A}_K$ with almost same efficiency as that of $\exp(t\hat{A}_L)$.
The natural domain of application of (\ref{nefedov:2.8}),
(\ref{nefedov:2.81}) is the one when $\hat{A}$ is a difference operator with
constant coefficients in $C({\mathbb Z}^s)$, $s\geq 1$ (i.e., the
functions $a_x(y)$ in expression (\ref{nefedov:2.1}) only depend on the
difference $x-y$), and $\Omega=\prod\limits_{j=1}^s\{0, 1, \ldots, N_j-1\}$ is a
parallelepiped in ${\mathbb Z}^s$.  In this case there is a specific
extension operator $\hat{L}$ which is defined by the relation
$(\hat{L}f)(x)=f(x\mod  N)$, where $(x\mod N)_j=x_j\mod N_j$ for $j=0,
1, \ldots, s$, corresponding to the periodic boundary conditions for
$\hat{A}$.  The exact value of operator $\exp(t\hat{A}_K)$ is
calculated using a multidimensional discrete Fourier transform, the
calculation procedure involves about $M\log_2M$ operations, where
$M=N_1\cdot\ldots\cdot N_s$.  The number of points of the set
$\partial_A\Omega$ can be estimated as
\[
|\partial_A\Omega|\leq C(A)M\sum_{j=1}^s\frac{1}{N_j},
\]
where the constant $C(A)$ depends on $\#(\mbox{supp}\;a(x))$.  Hence,
if $C(A)\ll\min N_j$, then the size of matrix $\hat{G}_{KL}$ is much
smaller than that of $\hat{A}_K$; this allows to effectively use
formulas (\ref{nefedov:2.8}), (\ref{nefedov:2.81}).  Such a situation
occurs in approximations of differential operators with difference
ones, and the constant $C(A)$ in this case depends on the order of the
operator approximated and, generally, on a method of approximation.

\section{Error and stability of the method}

Let us now study error and stability of numerical algorithms based on
formulas~(\ref{nefedov:2.8}) and~(\ref{nefedov:2.81}).

Let $f(t, \cdot)=\exp(t\hat{A}_K)g$ be the exact solution of problem
(\ref{nefedov:2.2}) with operator $\hat{A}_K$ and $h_j(t,
\cdot)=\hat{S}_j(t)g$, $j=1, 2$, $g\in C(\Omega)$.  As an error
estimate of one step of the OE algorithms we consider the norms of
differences of the functions $f(t, \cdot)$ and $h_j(t, \cdot)$:
\[
\varepsilon_j(t, g)=\|f(t, \cdot)-h_j(t, \cdot)\|, \qquad j=1, 2,
\qquad g\in C(\Omega).
\]

By expanding $f-h_j$ in the Taylor series at $t=0$ we
find
\begin{gather}
\varepsilon_1(t, g)=\left\|\left[\hat{A}_K, \hat{G}_{KL}\right]g\right\|
                 \frac{t^2}{2}+O\left(t^3\right),
\label{nefedov:e1}
\\
\varepsilon_2(t, g)=\left\|\left(
             \left[\hat{A}_K, \left[\hat{A}_K, \hat{G}_{KL}\right]\right]-
  \frac{1}{2}\left[\hat{G}_{KL}, \left[\hat{A}_K, \hat{G}_{KL}\right]\right]
                   \right)g\right\|
                 \frac{t^3}{12}+O\left(t^4\right),
\label{nefedov:e2}
\end{gather}
where $[\hat{A}, \hat{B}]=\hat{A}\hat{B}-\hat{B}\hat{A}$.

Thus, the schemes (\ref{nefedov:2.8}) and (\ref{nefedov:2.81}) have
the first and second order of approximation, respectively.

For comparison, consider similar estimates for the classical Euler and
Krank-Nickolson (KN) methods \cite{nefedov:1}.  The corresponding
approximations are of the form (cf.  (\ref{nefedov:2.8}) and
(\ref{nefedov:2.81}))
\begin{equation}\label{Euler}
\exp(t\hat{A}_K)=E+t\hat{A}_K+O(t^2)
\equiv \hat{S}_e (t)+O(t^2)
\end{equation}
for the Euler scheme and
\begin{equation}\label{KN}
\exp(t\hat{A}_K)=(E+\frac12t\hat{A}_K)(E-\frac12t\hat{A}_K)^{-1}+O(t^3) \equiv
\hat{S}_{kn} (t)+O(t^3)
\end{equation}
for the Krank-Nicolson one.  It is easy to obtain the estimates for
errors of these approximation, namely
\begin{gather}
\delta_1(t, g)=\|\hat{A}_K^2g\|\frac{t^2}{2}+O\left(t^3\right),
\label{nefedov:d1}
\\
\delta_2(t, g)=\|\hat{A}_K^3g\|\frac{t^3}{12}+O\left(t^4\right),
\label{nefedov:d2}
\end{gather}
for Euler and KN schemes, respectively.  Clearly, unlike the values
$\delta_j$, the estimates $\varepsilon_j$, where $j=1, 2$, are
determined by the norms of commutators of $\hat{A}_K$ with
$\hat{G}_{KL}$, rather than by the powers of $\hat{A}_K$.  This
accounts for the differences in the features of the OE algorithms and
classical schemes.

Stability analysis of the OE methods (\ref{nefedov:2.8}) and (\ref{nefedov:2.81})
requires evaluation of the norms of the relevant step operators
$\hat{S}_1(t)$ and $\hat{S}_2(t)$.  We are going to show that these
methods are stable for rather small $t>0$,  if operator $\hat{A}_K$
satisfies the condition

\begin{equation}\label{nefedov:Cond}
  \mbox{Re}(\hat{A}_Kg, g)<0 \qquad\text{for any } g\in C(\Omega).
\end{equation}
Indeed,  the functions
\[
s_j(t, g)=\|\hat{S}_jg\|^2,  \qquad g\in C(\Omega), \quad\text{for } j=1, 2
\]
are analytic in $t$ in a vicinity of zero, and $s_j(0, g)=\|g\|^2$.
Besides,
\[
\left.\frac{\partial\hat{S}_j}{\partial
t}\right|_{t=0}=2\,\mbox{Re}\,(\hat{A}_Kg, g)<0
\]
due to the assumption (\ref{nefedov:Cond}).
Therefore, for a sufficiently small positive values of $t$ the inequality
\[
s_j(t, g)=\|g\|^2+t\left.\frac{\partial \hat{S}_j}{\partial
t}\right|_{t=0} +O\left(t^2\right)\leq \|g\|^2,
\]
is valid, i.e., the schemes (\ref{nefedov:2.8}), (\ref{nefedov:2.81}) are
stable.  Note that the condition (\ref{nefedov:Cond}) means that the spectrum
$\mbox{spec}\; \hat{A}_K$ lies in the left half-plane~\cite{nefedov:9}; this
guarantees stability of the initial problem~(\ref{nefedov:2.2}).

There exist two classes of operators  $\hat{A}_K$
for which the OE method shows absolute stability.

1.  Schemes (\ref{nefedov:2.8}), (\ref{nefedov:2.81}) are absolutely
stable, if $\hat{A}_K$ and $\hat{A}_L$ are Hermitian operators and
$\mbox{spec}\;\hat{A}_L$ and $\mbox{spec}\;\hat{G}_{KL}$ both lie in
the left half-plane.  This immediately follows from the simplest
estimates:
\begin{gather*}
\|\hat{S}_1(t)\|\leq \|\exp(t\hat{G}_{KL})\|\|\exp(t\hat{A}_L)\|\leq 1,
\\
\|\hat{S}_2(t)\|\leq
\left\|\exp\left(\frac12t\hat{G}_{KL}\right)\right\|^2\|\exp(t\hat{A}_L)\|\leq 1.
\end{gather*}
This class of operators includes,  in particular,
the difference Laplace operator with the Dirichlet boundary conditions.

2.  The OE method is also absolutely stable,  if operators $\hat{A}_K$
and $\hat{A}_L$ are both skew-symmetric,  i.e.,  $\hat{A}_K=-\hat{A}_K^*$
and $\hat{A}_L=-\hat{A}_L^*$.  This condition is equivalent to the
case when $\text{spec}\hat{A}_K$ and $\text{spec}\hat{A}_L$ both lie on the
imaginary axis.  Then operators (\ref{nefedov:2.8}), (\ref{nefedov:2.81}) are unitary
(just as the operator $\exp(t\hat{A}_K)$); hence,  their norms are
equal to~1.  An example of such case is the Schr\"odinger operator.

The OE algorithm may, however, lack absolute stability even when both
$\mbox{spec}\; \hat{A}_K$ and $\mbox{spec}\; \hat{A}_L$ are in the
left half-plane.  This loss of stability is associated with the
positive eigenvalues available for the ``boundary" operator
$\hat{G}_{KL}$.  The example is a difference Laplacian with the
Neumann boundary conditions.

The above estimates can be illustrated by the $D$-expansion of Laplace
operator $\hat{\Delta}_D$ with the extension operator $\hat{D}$
corresponding to the Dirichlet boundary conditions ($\alpha=\beta=-1$
in (\ref{nefedov:2.7})).  In this case,  as mentioned earlier,  schemes
(\ref{nefedov:2.8}), (\ref{nefedov:2.81}) are absolutely stable.

The error estimates for the time step in the OE  methods (\ref{nefedov:e1}) and
(\ref{nefedov:e2}) depend on the initial vector $g$.  As typical
vectors $g$ we consider the eigenfunctions of the operator
$\hat{\Delta}_D$:
\[
\phi_j(k)=\sqrt{\frac{\sigma_j}{N}}
\sin\left(\frac{\pi(j+1)}{N}\left(k+\frac{1}{2}\right)\right), \qquad \mbox{for} \quad j, k=0, 1, \ldots, N-1,
\]
where $\sigma_j=2$,  $j=0, 1, \ldots, N-2$; $\sigma_{N-1}=1$. In this case
\begin{gather}
\varepsilon_1(t, \phi_j)\!=\!\left\{\!\!
   \begin{array}{lll}
   0 & \mbox{if}  &j=1, 3, \ldots, N-1, \\
  \displaystyle  \frac{4t^2}{2\sqrt{N}}
         \sqrt{-\sigma_j\mu_j(1+(3+\mu_j)^2)}+O\left(t^3\right)\!
         & \mbox{if}  &j=0, 2, \ldots, N-2,
   \end{array}
                        \right.
\label{nefedov:e3}\\
\varepsilon_2(t, \phi_j)\!=\!\left\{\!\!
   \begin{array}{lll}
   0 & \mbox{if} &j=1, 3, \ldots, N-1, \\
   \displaystyle\frac{t^3}{48\sqrt{N}}
\Bigl(-\sigma_j\mu_j
(1+2(3+\mu_j)^2\\
\qquad {}+(\mu_j^2+5\mu_j+7)^2)\Bigr)^{1/2}+O\left(t^4\right)
         & \mbox{if}  &j=0, 2, \ldots, N-2,
   \end{array}   \right.
\label{nefedov:e4}
\end{gather}
where $\mu_j=-4\sin^2\frac{\pi(j+1)}{2N}$ for $j=0, 1, \ldots, N-1$, are
the eigenvalues of the operator $\hat{\Delta}_D$.  The relevant
estimates for the Euler (\ref{Euler}) and KN (\ref{KN}) methods are of the form
\begin{gather}
\delta_1(t, g)=\mu_j^2\frac{t^2}{2}+O\left(t^3\right),
\label{nefedov:d3}
\\
\delta_2(t, g)=\mu_j^3\frac{t^3}{12}+O\left(t^4\right).
\label{nefedov:d4}
\end{gather}

Comparing the above estimates we see that, given the same order of
approximation, the error of the OE method is much greater for
eigenfunctions with small numbers and much smaller for higher
harmonics.

\resetfootnoterule

Observe that the error in the classical schemes comes from the
difference between the eigenvalues of the step operator and
$\exp(t\hat{A}_K)$, while their eigenfunctions coincide.  The step
operators of OE algorithm have error in both the eigenvalues and the
eigenfunctions\footnote{It is interesting to observe that the spectra
of $\hat{S}_1(t)$ and $\hat{S}_2(t)$ coincide and their eigenfunctions
differ only at boundary points.}.  However, as shown numerically, the
eigenvalues of the operators $\hat{S}_1(t)$ and $\hat{S}_2(t)$
approximate the spectrum of $\exp(t\hat{A}_K)$ better than the
eigenvalues of the classical schemes.

To make sure the above is true, let us find spectrum $\lambda_j(t)$, $j=0,
1, \ldots, N-1$, of the step operator $\hat{S}_2(t)$ for Laplacian
$\hat{\Delta}_D$.  We show in Appendix that for $j$ odd the
eigenvalues and eigenfunctions of operators $\hat{S}_2(t)$ and
$\exp(t\hat{\Delta}_D)$ coincide.  This is exactly why for these
harmonics the errors of the OE methods (\ref{nefedov:e3}) and (\ref{nefedov:e4})
vanish.  The remaining $\frac{N}{2}$ eigenvalues are
$\lambda_{2j}(t)=\exp(t\xi_j(t))$, $j=0, 1, \ldots, \frac{N}{2}-1$, where
$\xi_j(t)$ fulfills the ``dispersion" equation (\ref{nefedov:p10}).

The results of numerical solution of this equation are
shown in Fig.~1 which provides the values of $|\xi_j-\mu_{2j}|$ as
function of $j$ (curve~1).

Clearly,  for the majority of harmonics the eigenvalue error of the OE
method is much smaller than the error of the KN scheme (curve~2).
Besides,  one should note the ``uniformity" of the spectrum estimate of
scheme (\ref{nefedov:2.81}): the error weakly depends on the number of the
eigenvalues.  A similar situation holds also for the operator
$i\hat{\Delta}_D$.


\begin{figure}[!t]
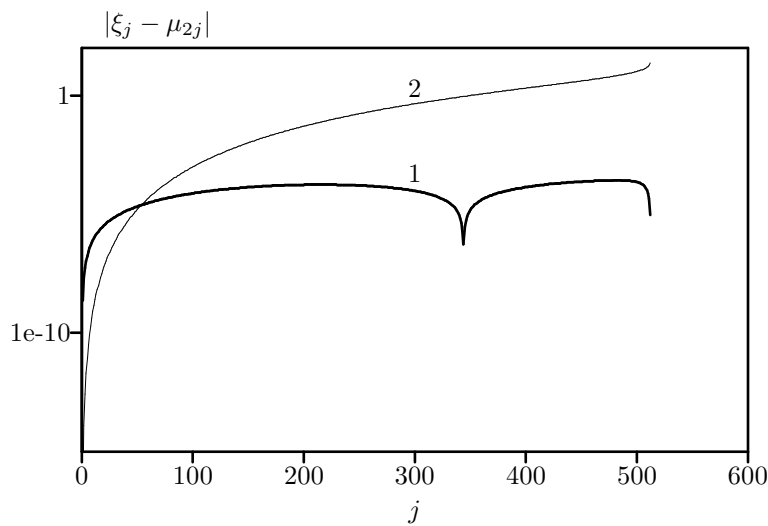

 \centering\include{fig1_nef}

\vspace{-5mm}

\caption{Error in even eigenvalues of step operator for
the OE method~(1) and the Krank--Nickolson method~(2)
for one-dimensional Laplace operator with
Dirichlet boundary conditions $t=\frac12$,  $N=1024$.}
\end{figure}

\begin{figure}[!t]
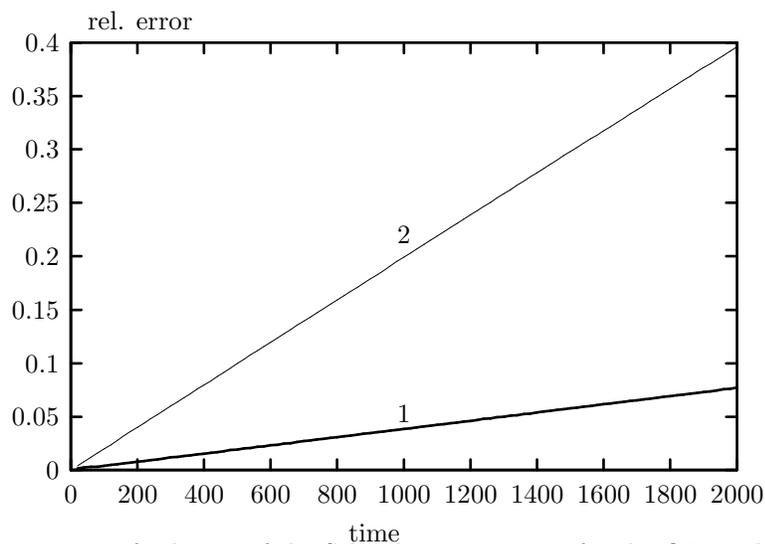



\centering\include{fig2_nef}

\vspace{-5mm}

\caption{Relative error of solution of the Schr\"odinger equation
for the OE method (1) and the Krank--Nickolson method (2)
with random initial vector,  $t=\frac12$,  $N=128$.}
\end{figure}

This property of the step operator in the numerical scheme is
important when solution of the input evolution problem includes
contributions from all eigenfunctions of operator~$\hat{A}_K$.  This
is the case, for example, in solving problem (\ref{nefedov:2.4}) with
skew-Hermitian operator~$\hat{A}_K$ (Schr\"odinger equation).

Fig.~2 shows relative errors of the numerical solution of the problem
\[
\frac{\partial f}{\partial t}=i\hat{\Delta}_D f, \qquad f(0)=g
\]
as a function of time for the OE scheme (\ref{nefedov:2.81}) (curve~1)
and the Krank--Nickolson one (curve~2).  The initial vector $g$ is
chosen as a random one; it is uniformly distributed on the unit sphere
in ${\mathbb C}^N$.

\section{Conclusion}

Splitting methods, including the operator exponential one, are widely
used for solving difference linear and quasilinear initial-boundary
value problems \cite{nefedov:2,nefedov:4}.  The proposed modification
of the OE method can be applied when the evolution operator is represented
as a sum of a difference operator with constant coefficients on a
rectangular domain in ${\mathbb Z}^s$ (the main part) and some,
perhaps nonlinear, operator (perturbation).  If the boundary
conditions for the main part do not allow explicit computation of the
input operator exponential, the problem can be approached by a
splitting method in two stages: first we split off the perturbation,
and then calculate an approximate exponential of the main part using
the method proposed in this work.

Consider a rectangular domain in ${\mathbb Z}^2$ for the Laplacian with
boundary conditions of the 3rd kind.  Even in this case application of
methods like the implicit Euler or Krank--Nickolson schemes requires
iteration procedures to obtain the resolvent.  The method we propose
is explicit and, as follows from the one-dimensional examples provided
in the work, competitive with conventional methods.

One of important features of our method is a ``uniform" property of
the spectral estimate of the initial problem.  We have succeeded in
applying the scheme described to solve one- and two-dimensional
Ginzburg--Landau equation~\cite{nefedov:10}, four-order diffusion
equation~\cite{nefedov:11} and some other problems.

\appendix
\section*{Appendix}
\renewcommand{\theequation}{{\rm A}.\arabic{equation}}

In what follows we derive the equation for the eigenvalues of operator
$\hat{S}_2(t)$ defined by relation (\ref{nefedov:2.81}) for one-dimensional
difference Laplace operator $\hat{\Delta}_D$ with the Dirichlet
boundary conditions.  To simplify the calculations, we assume that $N$
is even: $N=2M$.

In this case the ``boundary" operator $\hat{G}_{KL}$ is of the form:
$\hat{G}_{KL}=-2\hat{Q}_0, $ where $\hat{Q}_0$ is the orthogonal
projection on vector $x_0=\frac{e_1+e_{N-1}}{\sqrt{2}}$, where
$\{e_j\}_{j=0}^{N-1}$ is a standard basis in~${\mathbb C}^N$.  It is
easy to calculate the exponential of such an operator:
\begin{equation}
\exp\left(\frac12t\hat{G}_{KL}\right)=\hat{E}+\left(e^{-t}-1\right)\hat{Q}_0.
\label{nefedov:p2}
\end{equation}

Let $\lambda$ be an eigenvalue of $\hat{S}_2(t)$ and $\Phi_\lambda$
the corresponding eigenfunction.  Then
\begin{equation}
\hat{S}_2(t)\Phi_\lambda\equiv
\exp\left(\frac12t\hat{G}_{KL}\right)\exp\left(t\hat{\Delta}_D\right)\exp
\left(\frac12t\hat{G}_{KL}\right)
\Phi_\lambda=\lambda\Phi_\lambda.
\end{equation}
Using expression (\ref{nefedov:p2}) and notation
$\Psi_\lambda=\exp\left(\frac12t\hat{G}_{KL}\right)\Phi_\lambda$,
we express the relation in the form:
\begin{equation}
\left(\lambda\hat{E}-\exp(t\hat{\Delta}_D)\right)\Psi_\lambda=
\lambda\left(1-e^{2t}\right)\hat{Q}_0\Psi_\lambda.
\label{nefedov:p4}
\end{equation}
Let $\{f_j\}_{j=0}^{N-1}$ be the eigenbasis of operator
$\hat{\Delta}_L$ (Laplacian with periodic boundary conditions):
\begin{equation}
\left(f_j\right)(k)=\frac{1}{\sqrt{N}}
 \exp\left(i\frac{2\pi jk}{N}\right), \qquad j, k=0, 1, \ldots, N-1,
\label{nefedov:p5}
\end{equation}
the corresponding eigenvalues $\nu_j$ being
\begin{equation}
\nu_j=-4\sin^2\frac{\pi j}{N}, \qquad j=0, 1, \ldots, N-1.
\label{nefedov:p6}
\end{equation}
Note that $\nu_j=\nu_{N-j}$ and $\nu_j=\mu_{2j+1}$ for $j=1, 2,
\ldots, M$, where $\mu_k$ are the eigenvalues of $\hat{\Delta}_D$.
Note also that the right hand side of (\ref{nefedov:p4}) for
any value of $\Psi_\lambda$ is proportional to vector $x_0$ whose
expansion with respect to basis (\ref{nefedov:p5}) is:
\begin{equation}
x_0=\sum_{j=0}^{N-1}(x_0, f_j)f_j=
 \frac{1}{\sqrt{2N}}\sum_{j=0}^{N-1}
 \left(1+\exp\left(-i\frac{2\pi j}{N}\right)\right)f_j.
\label{nefedov:p7}
\end{equation}
If $c_j$, where $j=0, 1,\ldots, N-1$, are the expansion coefficients
of function $\Psi_\lambda$ with respect to basis (\ref{nefedov:p5}),
then relation (\ref{nefedov:p4}) can be expressed in the form:
\begin{equation}
(\lambda-\exp(t\nu_j))c_j=a_\lambda\lambda\left(1-e^{2t}\right)(x_0, f_j), \qquad j=0, 1, \ldots, N-1,
\label{nefedov:p8}
\end{equation}
where $a_\lambda=(\Psi_\lambda, x_0)$.

In order to find all solution of equation (\ref{nefedov:p7}) we
consider two cases:

1) $a_\lambda=0$.  In this case $c_j=0$ or $\lambda=\exp(t\nu_j)$ for
each $j=0, 1, \ldots, N-1$.  Since $\Psi_\lambda\neq 0$, there exists
index $l$ such that $c_l\neq 0$.  In this case
$\lambda=\exp(t\nu_l)=\exp(t\nu_{N-l})$ for $l=1, 2, \ldots, M$ and
$c_j=0$ for $j\neq l, N-l$.  The corresponding eigenfunctions are
found from the condition $a_\lambda=0$:
\begin{equation}
\Psi_l(k)\sim \sin\frac{2\pi l}{N}\left(k+\frac{1}{2}\right), \quad k=0, 1, \ldots, N-1,
\qquad l=1, 2, \ldots, M.
\end{equation}
Note that these functions are the eigenfunctions of operator
$\hat{\Delta}_D$ corresponding to the eigenvalues $\mu_{2l-1}$, where
$l=1, 2, \ldots, M$.

2) $a_\lambda\neq 0$.
In this case $\lambda\neq\exp(t\nu_j)$ for each $j=0, 1, \ldots, N-1$
as follows from (\ref{nefedov:p8}).
Multiplying both sides of (\ref{nefedov:p7}) by
$\frac{(f_j,  x_0)}{\lambda-\exp(t\nu_j)}$ and summing over $j$
we obtain:
\[
a_\lambda=a_\lambda\lambda\left(1-e^{2t}\right)\sum_{j=0}^{N-1}
   \frac{|(x_0, f)|^2}{\lambda-\exp(t\nu_j)}.
\]
With condition $a_\lambda\neq 0$ and also (\ref{nefedov:p6})
and (\ref{nefedov:p7}) this relation takes the form:
\begin{equation}
1=\frac{1-e^{2t}}{M}\left(\frac{1}{1-\exp(t\xi)}+\frac{1}{2}
  \sum_{j=1}^{M-1}\frac{4+\nu_j}{1-\exp(t(\nu_j-\xi))}\right),
\label{nefedov:p10}
\end{equation}
where $\xi=\ln\frac{\lambda}{t}$.  It is easy to see that equation
(\ref{nefedov:p10}) has exactly $M$ real roots, one in each interval $(\nu_j,
\nu_{j+1})$, where $j=0, 1, \ldots, M-1$. These roots can be easily found
numerically by the bisection method. 

Denote the solutions of (\ref{nefedov:p10}) by $\xi_j(t)$, where $j=1, 2, \ldots,
M$.  Then, for the spectrum $\lambda_j(t)$, where $j=0, 1, \ldots, N-1$,
of operator $\hat{S}_2(t)$ we finally obtain:
\[
\lambda_{2j}(t)=\exp(t\xi_j(t)), \qquad \lambda_{2j+1}=\exp(t\nu_j),
\qquad j=0, 1, \ldots, M-1.
\]

Note that for the eigenvalues $\lambda(t)=\exp(i t\xi(t))$ of the step
operator $\hat{S}_2(t)$ corresponding to (Schr\"odinger) operator
$i\hat{\Delta}_D$ equation (\ref{nefedov:p10}) is of the form:
\[
1=\frac{\tan
t}{M}\left(-\frac{1}{\tan\left(\frac12t\xi\right)}+\frac{1}{2}
\sum_{j=1}^{M-1}\frac{4+\nu_j}{\tan\left(\frac12t\left(\nu_j-\xi\right)\right)}\right).
\]
The solutions of this equation are in good agreement with those of
(\ref{nefedov:p10}).

\subsection*{Acknowledgements}

We are thankful to E~I~Gordon and M~A~Antonets for
fruitful discussions.  The work is partially supported by the RFBR
grant 99-02-16188.

\label{nefedov-lastpage}


\begin{thebibliography}{99}
\small

\bibitem{nefedov:1}
Samarsky A~A and Gulin  A~V,  Numerical Methods,  Moscow -  Nauka,
1989.
\bibitem{nefedov:2}
Hardin R~H and  Tappert F~D,  Application of the Split-Step Fourier
Method to the Numerical Solution of Nonlinear and Variable Coefficient
Wave Equations, {\it  SIAM Rev.  Chronicle}  {\bf 16} (1973), 423.

\bibitem{nefedov:3}
Kato T,  Perturbation Theory for Linear Operators,
Springer-Verlag -  Berlin-Heidelberg-New York,  1966.
\bibitem{nefedov:4}
Marchuk G~I,  Methods of Computational Mathematics,  Nauka -
Moscow,  1989.
\bibitem{5}
Nefedov I~M and Shereshevskii~I~A, On the Calculation of the Exponential of
Difference Operators, {\it Mat.  Model.} {\bf 7}, N~5 (1995), 88. (Russian)

\bibitem{nefedov:6}
Akhiezer N~I and  Glazman I~N,  Theory of Linear Operators in
Hilbert Space, Second revised and augmented edition, Nauka -  Moscow
1966 (in Russian); Third edition,  corrected and augmented.
Vshcha Shkola,  Kharkov,  Vol.  I,  1977,  Vol.  II,  1978 (in Russian);
Akhiezer  N~I and Glazman  I~M, Theory of
Linear Operators in Hilbert Space,  Translated from the Russian and
with a preface by Merlynd Nestell,  Reprint of the 1961 and 1963
translations,  Two volumes bound as one,  Dover Publications,  Inc. -
New York,  1993.
\bibitem{nefedov:7}
Albeverio S, Gesztesy F, Hegh-Krohn R. and Holden H,
Solvable Models in Quantum Mechanics,  Texts and Monographs in
Physics,  Springer-Verlag - New York-Berlin, 1988.

\bibitem{nefedov:8}
Shereshevskii~I~A,  A Finite Dimensional Analog of the Krein Formila,
submitted to JNMP.

\bibitem{nefedov:9}
Riesz F and Sz.-Nagy B,  Le\c{c}ons D'analyze fonctionnelle,
Akademiai Kiad\'o -  Budapest,  1972.
\bibitem{nefedov:10}
Andronov A,  Gordion I,  Kurin V,  Nefedov I and Shereshevsky I,
Kinematic Vortices and Phase Slip Lines in the Dynamics of the
Resistive State of Narrow Suprconductive thin Film Channels,
{\it Physica C}  {\bf 213} (1993), 193.
\bibitem{nefedov:11}
Fraerman A~A,  Mel'nikov~A~S,  Nefedov~I~M, Shereshevskii~I~A and
Shpiro A~V,  Nonlinear Relaxation Dynamics in Decomposing Alloys:
One-Dimensional Cahn--Hillard Model,  {\it Phys. Rev.}  {\bf 10}, N~10 (1997),  6316--6323.

\end{thebibliography}
\end{document}